\documentclass{amsart}
\usepackage{amscd,amssymb}
\usepackage{fullpage}
\usepackage{righttag}
\usepackage{pb-diagram,lamsarrow,pb-lams}

\newcommand{\ov}[1]     {\overline{#1}}

\newcommand{\gbin}[2]
 {\left[\begin{array}{c}{#1}\\{#2}\end{array}\right]_p}
\newcommand{\bsm}       {\left(\begin{smallmatrix}}
\newcommand{\esm}       {\end{smallmatrix}\right)}

\newcommand{\ds}        {DS^0}
\newcommand{\bs}[1]     {B\Sigma_{#1}}
\newcommand{\bsp}[1]    {B\Sigma_{#1 +}}
\newcommand{\blob}      {\bullet}
\newcommand{\eds}       {E^0\ds}
\newcommand{\ev}        {E^\vee}
\newcommand{\evds}      {E^\vee\ds}
\newcommand{\kds}       {K^0\ds}
\newcommand{\kvds}      {K_0\ds}
\newcommand{\psis}      {\psi_*}
\newcommand{\psit}      {\psi_\times}
\newcommand{\psib}      {\psi_\bullet}
\newcommand{\mxi}       {\mathfrak{m}}
\newcommand{\Gm}        {\Gamma}
\newcommand{\Lm}        {\Lambda}
\newcommand{\Oi}        {\Omega^\infty}
\newcommand{\Sgi}       {\Sigma^\infty}

\newcommand{\oibp}      {\Omega^\infty BP}
\newcommand{\db}        {\overline{d}}
\newcommand{\un}[1]     {\underline{#1}}

\newcommand{\tm}        {\times}
\newcommand{\tX}        {\times_X}
\newcommand{\ot}        {\otimes}

\newcommand{\ote}       {\otimes_E}

\newcommand{\hote}      {\widehat{\otimes}_E}
\newcommand{\Wedge}     {\vee}
\newcommand{\Smash}     {\wedge}
\newcommand{\Lhat}      {{\widehat{L}}}

\newcommand{\CC}        {\mathcal{C}}
\newcommand{\BD}        {{\mathbb{D}}}

\newcommand{\BC}        {B\mathcal{C}}
\newcommand{\cpi}       {\mathbb{C} P^\infty}
\newcommand{\OO}        {\mathcal{O}}

\newcommand{\Level}     {\operatorname{Level}}
\newcommand{\Div}       {\operatorname{Div}}
\newcommand{\Hom}       {\operatorname{Hom}}

\newcommand{\Ind}       {\operatorname{Ind}}
\newcommand{\Prim}      {\operatorname{Prim}}
\newcommand{\Sub}       {\operatorname{Sub}}
\newcommand{\Tor}       {\operatorname{Tor}}
\newcommand{\ann}       {\operatorname{ann}}
\newcommand{\image}     {\operatorname{image}}
\newcommand{\res}       {\operatorname{res}}
\newcommand{\tr}        {\operatorname{tr}}
\newcommand{\spf}       {\operatorname{spf}}

\newcommand{\st}        {\;|\;}
\newcommand{\fpl}       {[\![} 
\newcommand{\fpr}       {]\!]} 
\newcommand{\fps}[2]    {{#1 \fpl #2 \fpr}} 
\newcommand{\sse}       {\subseteq}

\newcommand{\burn}      {\mathbb{B}}          
\newcommand{\latt}      {\mathbb{L}}          
\newcommand{\real}      {\mathbb{R}}              
\newcommand{\rinf}      {\mathbb{R}^\infty}
\newcommand{\cplx}      {\mathbb{C}}     
\newcommand{\QZp}       {{\mathbb{Q}_p/\mathbb{Z}_p}}         
\newcommand{\Zh}        {\mathbb{Z}}              
\newcommand{\Zp}        {\mathbb{Z}_p}   
\newcommand{\Qp}        {\mathbb{Q}_p}         
\newcommand{\Fp}        {{\mathbb{F}_p}}          
\newcommand{\Fq}        {{\mathbb{F}_{p^n}}}
\newcommand{\WF}        {W\mathbb{F}_{p^n}}              
\newcommand{\GG}        {\mathbb{G}}

\newcommand{\al}        {\alpha}
\newcommand{\bt}        {\beta} 
\newcommand{\dl}        {\delta}
\newcommand{\ep}        {\epsilon}
\newcommand{\sg}        {\sigma}
\newcommand{\tht}       {\theta}
\newcommand{\Dl}        {\Delta}
\newcommand{\Sg}        {\Sigma}
\newcommand{\Om}        {\Omega}

\newcommand{\xra}       {\xrightarrow}

\newcommand{\era}       {\twoheadrightarrow}
\newcommand{\mra}       {\rightarrowtail}

\newcommand{\colim}{\operatornamewithlimits{\underset{\longrightarrow}{lim}}}

\renewcommand{\:}{\colon}

\newtheorem{theorem}{Theorem}[section]

\newtheorem{lemma}[theorem]{Lemma}
\newtheorem{proposition}[theorem]{Proposition}
\newtheorem{corollary}[theorem]{Corollary}
\theoremstyle{definition}
\newtheorem{remark}[theorem]{Remark}
\newtheorem{definition}[theorem]{Definition}

\newenvironment{diag}{
 \renewcommand{\typeout}[1]{}
 \begin{displaymath}
 \begin{diagram}}{
 \end{diagram}
 \end{displaymath}}

\begin{document}
\title{Morava {$E$}-Theory of Symmetric Groups}
\author{N.~P.~Strickland}
\address{Trinity College, Cambridge CB2 1TQ, England}
\email{neil@dpmms.cam.ac.uk}
\bibliographystyle{abbrv}

\begin{abstract}
 We compute the completed $E(n)$ cohomology of the classifying spaces
 of the symmetric groups, and relate the answer to the theory of
 finite subgroups of formal groups.
\end{abstract}
\maketitle 

\section{Introduction}\label{sec-intro}

Fix a prime $p$ and an integer $n>0$.  In this paper we study the
$\widehat{E(n)}$-cohomology of the space
\[ \ds = \bigvee_{k\ge 0} \bs{k+}, \]
where $\Sg_k$ denotes the symmetric group on $k$ letters.  It will be
technically convenient to work with a certain $2$-periodic algebraic
extension of $\widehat{E(n)}$ called Morava $E$-theory, which we
denote by $E$.  See section~\ref{sec-mke} for details.

It will turn out that $\eds$ is related to the theory of subgroups of
formal groups.  This general theme has its roots in work of Mike
Hopkins and Matthew Ando~\cite{an:ifg}.  We will also make heavy use
of results of Takuji Kashiwabara~\cite{ka:bpo}, and the generalised
character theory of Hopkins, Kuhn and Ravenel~\cite{hokura:ggc}.

Hopkins has proved that $E$ can be made into an $E_\infty$ ring
spectrum, but unfortunately no proof has yet appeared.  Our
calculation of $\eds$ is an important first step towards the
exploitation of the $E_\infty$ structure, as will be explained in a
sequel to this paper.  We have avoided using the $E_\infty$ structure
here, because of the lack of a published account.

\subsection{General notation and conventions}
\label{sec-not}

We write $\Fq$ for the finite field with $p^n$ elements; there is a
unique such field up to unnatural isomorphism.  In particular
$\Fp=\Zh/p$.  We write $\Zp$ for the ring of $p$-adic integers, and
$\Qp$ for its field of fractions, so $\Qp/\Zp\simeq\colim_n\Zh/p^n$.

Unless otherwise stated, all our vector bundles, group representations
and tensor products are over the complex numbers.  Given a
representation $V$ of a finite group $G$, we use the same symbol $V$
to denote the associated vector bundle $EG\tm_GV$ over $BG$.  Given a
vector bundle $V$ over a space $X$, we write $X^V$ for the associated
Thom space, or for its suspension spectrum.  Given an integer $n$, we
sometimes write $n$ for the trivial bundle of complex dimension $n$,
so that $n\ot V$ is the direct sum of $n$ copies of $V$.  

Because $E^*$ is periodic, we can choose a degree-zero complex
orientation $x\in\tilde{E}^0\cpi$.  Using this, we define Thom
classes, Euler classes and Chern classes for complex bundles, in the
usual way (all of them in degree zero).  We write $e(V)$ for the Euler
class of a bundle $V$.

We write $V_m$ for the usual permutation representation of $\Sg_m$ on
$\cplx^m$, or for the corresponding bundle $E\Sg_m\tm_{\Sg_m}\cplx^m$
over $B\Sg_m$.  This contains an evident copy of the trivial
representation of complex dimension one, so we can interpret $V_m-1$
as an honest complex representation or bundle, and thus define the
Euler class $c_m=e(V_m-1)\in E^0\bs{m}$.

\subsection{Statement of results}
\label{sec-results}

The stable transfer maps $\bs{k+l}\xra{}\bs{k}\tm\bs{l}$ and the
diagonal maps $\bs{k}\xra{}\bs{k}\tm\bs{k}$ give two products on
$\eds$, written $\tm$ and $\blob$.  The evident maps
$\bs{k}\tm\bs{l}\xra{}\bs{k+l}$ give rise to a coproduct $\psis$ on
$\eds$ (provided that we use a completed tensor product).

The following theorem summarises most of our results.
\begin{theorem}
 With the above structure, $\eds$ is a Hopf ring.  It is a formal
 power series ring under $\tm$.  The $\tm$-indecomposables are
 \[ \Ind(\eds) = \prod_{k\ge 0} \ov{R}_k, \]
 where 
 \[ \ov{R}_k =
     E^0\bs{p^k}/\operatorname{transfer}(E^0B\Sg_{p^{k-1}}^p).
 \]
 This ring is naturally isomorphic to the ring $\OO_{\Sub_k(\GG)}$
 studied in~\cite{st:fsf}, which classifies subgroup-schemes of degree
 $p^k$ in the formal group $\GG$ associated to $E^0\cpi$.  It is a
 Gorenstein local ring, and a free module over $E^0$, with rank given
 by the Gaussian binomial coefficient
 \[ \gbin{n+k-1}{n-1} = \prod_{j=1}^{n-1} \frac{p^{k+j}-1}{p^j-1}
    = |\{ \text{ subgroups } A < (\QZp)^n \text{ with } |A|=p^k\}|.
 \]
 In Section~\ref{sec-basis}, we will give an explicit basis consisting
 of monomials in the Chern classes of $V_{p^k}$.  We also have
 \[ \Prim(\eds) =
    \prod_k\ker(\res\:E^0\bs{p^k} \xra{} E^0B\Sg_{p^{k-1}}^p).
 \]
 This is a free module over $\Ind(\eds)$ on one generator $c$, whose
 component in $E^0\bs{m}$ is the Euler class $c_m=e(V_m-1)$.
\end{theorem}

We now give a brief outline of the methods used to prove this.
Firstly, we appeal to the results of~\cite{sttu:rme}.  There, we
showed that the maps of spectra that we used to define the products
and coproduct on $\eds$ interact in the right way to make $\eds$ a
Hopf ring; this is explained in Section~\ref{sec-epf}.  Also,
following work of Hopkins, Kuhn and Ravenel~\cite{hokura:ggc}, we
defined an $E^0$-algebra $L$ which is a free module over $p^{-1}E^0$,
and we gave a complete description of $L\ote E^0\bs{k}$.  This
essentially tells us everything we need to know about
$p^{-1}E^0\bs{k}$.  It is also known that $K^{\text{odd}}\bs{k}=0$
(see~\cite[Section 6]{hokura:ggc} or~\cite{hu:mkw}).  We recall the
proof in Section~\ref{sec-EBSk}, because we need to extract from it
the fact that $c_{p^k}^{(p^n-1)/(p-1)}\neq 0$ in $K^0\bs{p^k}$, which
turns out to be crucial later on.  It is not hard to deduce from the
fact that $K^{\text{odd}}\bs{k}=0$ that $E^*\bs{k}$ is a finitely
generated free module over $E^*$, concentrated in even degrees.  This
of course means that the loss of information in passing from
$E^0\bs{k}$ to $p^{-1}E^0\bs{k}$ is relatively small.

In Section~\ref{sec-kvds} we recall Kashiwabara's
result~\cite[Corollary 6.4]{ka:bpo} that $\kvds$ is a polynomial ring
under the product dual to $\psis$, and we give a simplified proof.  It
follows that $\evds$ is a completed polynomial ring over $E^0$, and
thus that $\Ind(\evds)$ is a completed free module and is a summand in
$\evds$.  Dually, we conclude that $\Prim(\eds)$ is a product of free
modules and is a summand in $\eds$.

We next consider $E^0QS^2$.  There is only one sensible product on
$E^0QS^2$, which comes from the diagonal map $QS^2\xra{}QS^2\tm QS^2$,
which is also the map obtained by applying $Q$ to the pinch map
$S^2\xra{}S^2\Wedge S^2$.  We also prove in Section~\ref{sec-kvds}
that $E^0QS^2$ is a formal power series ring under this product, using
an unpublished argument of Kashiwabara.

Our next task is to transfer this result to $\eds$, using the Snaith
splitting $\Sgi QS^2=\Sgi DS^2$ and the Thom isomorphism
$\tilde{E}^0D_k(S^2)=\tilde{E}^0\bs{k}^{V_k}=E^0\bs{k}$.  We need a
property of the Snaith splitting that does not seem to appear in the
literature; in order to prove it, we recall the details of the
construction in Section~\ref{sec-snaith}.  We find that our
isomorphism converts the product on $E^0QS^2$ to the $\tm$-product on
$\eds$, so we conclude that $\eds$ is a formal power series ring
under the $\tm$-product.  It follows from this that $\Ind(\eds)$ is a
product of free modules over $E^0$, and a retract of $\eds$.  

In Section~\ref{sec-prim-ind}, we study the indecomposables and
primitives more closely.  We distinguish between two different kinds
of primitives, and we define indecomposables and primitives in more
elementary terms.  We find that the double suspension map
$E^*QS^2\xra{}E^*QS^0$ corresponds to multiplication by $c_k=e(V_k-1)$
in $E^0\bs{k}$.  This implies that $\Prim(\eds)$ is free of rank one
over $\Ind(\eds)$.

In Section~\ref{sec-subgp}, we use generalised character theory to
show that the representation $V_{p^k}$ defines a subgroup scheme of
the formal group $\GG$ over the scheme $\spf(\ov{R}_k)$.  This gives
us a classifying map $\OO_{\Sub_k(\GG)}\xra{}\ov{R}_k$, and we prove
that this is an isomorphism.  All this makes heavy use
of~\cite{st:fsf}.  

In Section~\ref{sec-basis}, we give a basis for $\OO_{\Sub_k(\GG)}$,
consisting of monomials in the Chern classes of $V_{p^k}$.  This is a
piece of pure algebra that should really have appeared
in~\cite{st:fsf}.  

\section{Morava $K$-theory and $E$-theory}
\label{sec-mke}

In this section, we give more precise definitions of the (co)homology
theories which we need.

We will consider a completed and extended version of $E(n)$, which we
call Morava $E$-theory.  We start with the graded ring
\[ E_* = \fps{\WF}{u_1,\ldots,u_{n-1}}[u^{\pm 1}] \]
(where $u_i\in E_0$ and $u\in E_2$).  We take $u_0=p$, $u_n=1$ and
$u_k=0$ for $k>n$.  We make $E_*$ into an algebra over $BP_*$ by
sending the Hazewinkel generator $v_k$ to $u^{p^k-1}u_k$.  This makes
$E_*$ into a Landweber exact $BP_*$-module, so we have a homology
theory
\[ E_*X = E_*\ot_{BP_*} BP_*X = E_*\ot_{MU_*} MU_*X. \]
This is represented by a spectrum $E$.  

There is also a ring spectrum $K$, with
\[ K_* = E_*/(u_0,\ldots,u_{n-1}) = \Fq[u^{\pm 1}] =
      \Fq\ot K(n)_*[u]/(u^{p^n-1}-v_n). 
\]
We shall call this the {\em Morava $K$-theory} spectrum.  There is a
ring map $E\xra{}K$ with the obvious effect on homotopy groups.  If
$p>2$ then $K$ is commutative.  If $p=2$ then there is a well-known
formula
\[ ab - ba = u Q_{n-1}(a) Q_{n-1}(b), \]
where $Q_{n-1}\:K\xra{}\Sg K$ is a derivation.  (A proof is given for
$K(n)$ in~\cite[Theorem 2.13]{st:pmm}, for example, and the statement
for $K$ follows easily.)  It follows that $K^0X$ is commutative
whenever $X$ is a space with $K^1X=0$.  This will be the case for all
spaces that we consider.

We also write $\Lhat X$ for the Bousfield localisation of a spectrum
$X$ with respect to $K(n)$ (or equivalently, with respect to $K$).  We
write
\[ \ev(X) = \pi_0 \Lhat(E\wedge X). \]
This is a more natural thing to consider than $E_0X$ when $X$ is
$K$-local.  In that case $E_0X$ is analogous to something like
$\Zp\ot\Zp$, which contains a large divisible subgroup, whereas 
$\ev X$ is analogous to the $p$-completion of $\Zp\ot\Zp$, which is
just $\Zp$.  It would lead us too far afield to justify these remarks
here; we refer the interested reader to~\cite{host:mkl} for some
discussion of related points.

\section{$E$-Theory of symmetric groups}\label{sec-EBSk}

In this section, we approach the calculation of $E^*\bs{k}$ using
fairly traditional methods.  We begin with a definition.

\begin{definition}\label{defn-Lambda}
 Define $\Lm=(\QZp)^n$, and $\Lm^*=\Hom(\Lm,\QZp)=\Zp^n$.  Let
 $\burn_k$ be the set of isomorphism classes of sets of order $k$ with
 an action of $\Lm^*$, and write $d(k)=|\burn_k|$.  Let $\latt_k$ be
 the set of lattices of index $p^k$ in $\Lm^*$, which bijects with the
 set of transitive $\Lm^*$-sets of order $p^k$, or with the set of
 subgroups of $\Lm$ of order $p^k$.  Write $\db(k)=|\latt_k|$.  It is
 shown in~\cite[Section 10.1]{st:fsf} that
 \[ \db(k) = \gbin{n+k-1}{n-1} = 
    \prod_{j=1}^{n-1} \frac{p^{k+j}-1}{p^j-1}.
 \]
\end{definition}

Our conclusions are as follows.
\begin{theorem}\label{thm-EBSk}
 $E^0\bs{k}$ is a Noetherian local ring and a free module over $E^0$
 of rank $d(k)$.  Moreover, we have
 \begin{align*}
  \ev\bs{k}     &= \Hom_{E^0}(E^0\bs{k},E^0)            \\
  K^0\bs{k}     &= \Fq\ote E^0\bs{k}                    \\
  K_0\bs{k}     &= \Hom_\Fq(K^0\bs{k},\Fq)              \\
  E^1\bs{k}     &= \ev_1\bs{k} = K^1\bs{k} = K_1\bs{k} = 0.
 \end{align*}
 Finally, if $k=p^m$ then $c_k^{(p^n-1)/(p-1)}$ is nonzero in
 $K^0\bs{k}$ (this will be important later).
\end{theorem}
\begin{proof}
 Assemble Proposition~\ref{prop-KWm}, Proposition~\ref{prop-ESk}, and
 Corollary~\ref{cor-tate}.  
\end{proof}

\begin{definition}\label{defn-Wm}
 Let $W_m$ be the $m$-fold iterated wreath product 
 $C_p\wr\ldots\wr C_p$.  It is well-known that this is the Sylow
 $p$-subgroup of $\Sg_{p^m}$, and that for any $k$ the Sylow subgroup
 of $\Sg_k$ is a product of $W_m$'s.
\end{definition}

\begin{proposition}\label{prop-KWm}
 The group $K^*BW_m$ is generated by transfers of Euler classes of
 complex representations of subgroups, and thus is concentrated in
 even degrees.  Moreover, the element $b_m=c_{p^m}^{(p^n-1)/(p-1)}$ is
 nonzero.
\end{proposition}
\begin{proof}
 The first sentence is proved in~\cite{hokura:ggc}; we need to recall
 the details in order to prove the second sentence.  

 Suppose that the proposition holds for $W_m$, this being trivial for
 $m=0$.  We can then choose a basis $\{e_1,\ldots,e_d\}$ for $K^0BW_m$
 over $\Fq$, where $e_1=b_m$ and each $e_k$ is a transferred Euler
 class.  There is a spectral sequence
 \[ H^*(C_p;K^*(BW_m)^{\ot p}) \Rightarrow K^*(BW_{m+1}). \]
 Let $z$ be a generator of $H^1(C_p;\Fp)$. 
 Write $y=\bt z\in H^2(C_p;\Fp)$, so that $y$ is the Euler class of
 the representation $L$ of $W_{m+1}$ associated to the character
 $W_{m+1}\xra{}C_p\subset S^1$.  Let $x\in K^0BW_{m+1}$ be the Morava
 $K$-theory Euler class for the same representation.

 Given a sequence of indices $I=(i_0\le\ldots\le i_{p-1})$ with
 $i_0<i_{p-1}$, let $e_I$ be the sum of the orbit under $C_p$ of
 $e_{i_0}\ot\ldots\ot e_{i_{p-1}}$.  We also write
 $e'_i=e_i\ot\ldots\ot e_i$.  The $E_2$ term is
 \[ E[z]\ot P[y]\ot \Fp\{e_I,e'_i\}/(x e_I,y e_I).\]
 We shall exhibit elements $\tilde{e}_I$ and $\tilde{e}'_i$ in
 $K^0BW_{m+1}$ that hit $e_I$ and $e'_i$ under the restriction map
 $K^*BW_{m+1}\xra{}BW_m^p$, which is the edge map of the spectral
 sequence.  First, we take
 \[ \tilde{e}_I = 
    \tr_{W_m^p}^{W_{m+1}}
     (e_{i_0}\ot\ldots\ot e_{i_{p-1}})\in K^0BW_{m+1}.
 \]
 It is easy to check that this hits $e_I$, and also that
 $x\tilde{e}_I=0=y\tilde{e}_I$.  Next, suppose that
 $e_i=\tr_H^{W_m}(e(V))$ for some complex representation $V$ of some
 subgroup $H\le W_m$.  There is an obvious way to extend the action of
 $H^p$ on $p\ot V$ to get an action of $C_p\wr H$.  Write
 \[ \tilde{e}'_i = \tr_{C_p\wr H}^{W_{m+1}}(e(p\ot V)) \in K^0(BW). \]
 This hits $e'_i$ under the edge map.

 Our spectral sequence is a module over the Atiyah-Hirzebruch spectral
 sequence  
 \[ H^*(BC_p;K^*)=E[z]\ot P[y]\Rightarrow K^*BC_p = P[x]/x^{p^n}. \]
 This has only one differential, viz. $d_{p^n-1}(z)=y^{p^n}$.  We
 conclude that $d_{p^n-1}(z e'_i)=y^{p^n}e'_i$, and after this
 differential the spectral sequence is concentrated in even bidegrees
 so nothing else can happen.  Thus, we have
 \[ K^*BW_{m+1} = K^*\{\tilde{e}_I,x^j \tilde{e}'_i\}, \]
 where $j$ runs from $0$ to $p^n-1$.  

 We still need to show that $b_{m+1}$ is nonzero.  In fact, we will
 show that it is a unit multiple of $x^{p^n-1}\tilde{e}'_1$.

 Let $U$ be the pull back of the regular representation of $C_p$ along
 the projection $W_{m+1}\era C_p$.  Note that
 $U\simeq\bigoplus_{r=0}^{p-1}L^{\ot r}$.  We have
 $e(L^{\ot r})=[r](x)$, where $[r](t)$ is the $r$-series for the
 formal group law associated to $E$.  This has the form
 $[r](t)=rt\pmod{t^2}$, so when $0<r<p$ we see that $[r](x)$ is a unit
 multiple of $x$.  It follows that $x^{p-1}$ is a unit multiple of
 $e(U-1)$, and thus that $x^{p^n-1}$ is a unit multiple of
 $e(U-1)^N=e(N\ot(U-1))$, where $N=(p^n-1)/(p-1)$.

 Next, recall that $b_m=e(V_{p^m}-1)^N=e(N\ot V_{p^m}-N)$.  The
 restriction of the representation $N\ot V_{p^{m+1}} - N\ot U$ of
 $W_{m+1}$ to the subgroup $W_m^p$ is isomorphic to
 $p\ot(N\ot V_{p^m} - N)$, so $\tilde{e}'_1=e(V_{p^{m+1}}-U)^N$.  It
 follows that $x^{p^n-1}\tilde{e}'_1$ is a unit multiple of
 $e(V_{p^{m+1}}-1)^N=b_{m+1}$, as required.
\end{proof}

\begin{proposition}\label{prop-filt}
 If $Z$ is a spectrum such that $E^*Z$ is finitely generated over
 $E^*$ and $K^*Z$ is concentrated in even degrees then $E^*Z$ is free
 and concentrated in even degrees.
\end{proposition}
\begin{proof}
 One can construct ring spectra $E/I_k$ (for $k>0$) with 
 \[ \pi_*(E/I_k)=E_*/I_k = \fps{\Fq}{u_k,\ldots,u_{n-1}}[u^{\pm 1}],\] 
 and cofibrations $E/I_k\xra{u_k}E/I_k\xra{}E/I_{k+1}$.  In
 particular, we have $E/I_n=K$, and we put $E/I_0=E$.  The
 cofibrations give short exact sequences
 \[ (E/I_k)^*(X)/u_k \mra (E/I_{k+1})^*(X) \era
     \ann(u_k,(E/I_k)^*(\Sg^{-1}X)).
 \]
 By induction on $k$, using the fact that $E^*$ is Noetherian, we see
 that $(E/I_k)^*X$ is finitely generated for all $k$.  We claim that
 it is actually free, and concentrated in even degrees.  This is
 trivially true for $(E/I_n)^*X$.  If it is true for $E/I_{k+1}$, we
 see from the short exact sequence that
 $(E/I_k)^{\text{odd}}(X)/u_k=0$, so $(E/I_k)^{\text{odd}}(X)=0$ by
 Nakayama's lemma.  We also see that
 $\ann(u_k,(E/I_k)^{\text{even}}(X))=0$, so $u_k$ acts regularly on
 $(E/I_k)^*X$.  Finally, we see that
 $(E/I_{k+1})^*X=(E/I_k)^*(X)/u_k$.  It follows by elementary
 commutative algebra that $(E/I_k)^*(X)$ is free, as required.
\end{proof}

\begin{proposition}\label{prop-ESk}
 $E^*\bs{k}$ is a free module of rank $d(k)$ over $E^*$, concentrated
 in even degrees.
\end{proposition}
\begin{proof}
 Because the Sylow subgroup of $\bs{k}$ is a product of $W_m$'s, we
 know from Proposition~\ref{prop-KWm} and a transfer argument that
 $K^*\bs{k}$ is concentrated in even degrees.  We also know
 from~\cite[Corollary 5.4]{grst:vlc} (for example) that $E^*BG$ is
 finitely generated as an $E^*$-module, for any finite group $G$.  It
 follows using Proposition~\ref{prop-filt} that $E^*\bs{k}$ is free
 and concentrated in even degrees.  The rank is necessarily the same
 as the rank of $p^{-1}E^*\bs{k}$ over $p^{-1}E^*$, and the main
 theorem of~\cite{hokura:ggc} shows that this is the same as the
 number of conjugacy classes of homomorphisms
 $\Lm^*=\Zp^n\xra{}\bs{k}$.  (A version of this theorem is proved
 in~\cite[Appendix A]{grst:vlc}, which may be easier to get hold of
 than~\cite{hokura:ggc}.)  It is easily seen that these biject with
 isomorphism classes of $\Lm^*$-sets of order $k$, so the rank is
 $d(k)$.  This is explained in more detail in~\cite{sttu:rme}.
\end{proof}

It turns out that whenever $E^*Z$ is a finitely generated free module
over $E^*$, the module $\ev_*Z$ is simply the dual of $E^*Z$.
However, we shall avoid having to prove this by giving a special
result for classifying spaces.
\begin{proposition}\label{prop-tate}
 If $G$ is a finite group then there is a weak equivalence
 \[ \Lhat(E\Smash BG_+)=F(BG_+,E), \]
 and thus an isomorphism $E^*BG=\ev_*BG$.
\end{proposition}
\begin{proof}
 Let $X$ be a spectrum.  The Greenlees-May theory of Tate
 spectra~\cite{grma:gtc} gives a cofibration of $G$-spectra
 \[ EG_+ \Smash X \xra{} F(EG_+,X) \xra{} t_G(X) . \]
 We now apply the Lewis-May fixed point functor and use Adams'
 isomorphism $(EG_+\Smash X)^G=EG_+\Smash_G X$ (see~\cite[Theorem
 7.1]{lemast:esh}). This gives a nonequivariant cofibration
 \[ BG_+ \Smash X \xra{} F(BG_+,X) \xra{} P_G(X) = t_G(X)^G, \]
 and thus a cofibration
 \[ \Lhat(BG_+\Smash X) \xra{} \Lhat F(BG_+,X) \xra{} \Lhat P_G(X). \]
 Moreover, all three terms here are exact functors of $X$.  Let $Y$ be
 a generalised Moore spectrum of type $n$, so that $E\Smash Y$ lies in
 the thick subcategory generated by $K$.  It is shown
 in~\cite{grsa:tcv,grsa:tct} that $P_GK=0$, so we must have
 $Y\Smash P_GE=0$.  As $K_*Y\neq 0$ we conclude that $K_*P_GE=0$ and
 thus $\Lhat P_GE=0$.  Note also that $E$ is $K$-local, and thus the
 same is true of $F(BG_+,E)$.  We thus conclude that
 $\Lhat(BG_+\Smash E)=F(BG_+,E)$, as claimed.
\end{proof}

\begin{corollary}\label{cor-tate}
 If $E^*BG$ is free then $\ev_*BG$ is its dual, and
 $K^*BG=K^*\ot_{E^*}E^*BG$.  
\end{corollary}
\begin{proof}
 Proposition~\ref{prop-tate} tells us that $\ev_*BG$ is free over
 $E^*$.  The cofibrations in the proof of Proposition~\ref{prop-filt}
 show us that $K^*BG=K^*\ot_{E^*}E^*BG$ and
 $K_*BG=K_*\ot_{E_*}\ev_*BG$.  Both of these are finitely generated
 over the graded field $K^*$, so they are dual to each other.  Thus,
 the natural map $d_K\:K_*BG\xra{}\Hom_{K^*}(K^*BG,E^*)$ is an
 isomorphism.  We also have a natural map
 $d_E\:\ev_*BG\xra{}\Hom_{E^*}(E^*BG,E^*)$, and one can check that
 $d_K=1_{K_*}\ot_{E_*}d_E$.  As $d_K$ is an isomorphism and the source
 and target of $d_E$ are free of finite rank over the local ring
 $E_*$, we conclude that $d_E$ is an isomorphism as required.
\end{proof}

\section{The extended power functor}\label{sec-epf}

We next recall some of the ideas developed in~\cite{sttu:rme}.  For
any based space $X$, we can define the total extended power
\[ DX = \bigvee_{k\ge 0} E\Sg_{k+}\Smash_{\Sg_k} X^{(k)}.
\]
In~\cite{lemast:esh}, it is shown how to define $DX$ when $X$ is a
spectrum, so that $D(\Sgi X)=\Sgi DX$.  We are mainly interested in
$\ds=\bigvee_{k\geq 0}\bsp{k}$.  We will sometimes think of this as a
space, and sometimes as a spectrum.  It can also be described as
$\Sgi\BC_+$, where $\CC$ is the category of finite sets.  We observed
in~\cite{sttu:rme} that it admits two products and two coproducts:
\[
\begin{array}{rl}
 \sg  \: & \ds \Smash \ds  \xra{}  \ds \\
 \mu  \: & \ds \Smash \ds  \xra{}  \ds \\
 \tht \: & \ds  \xra{}  \ds \Smash \ds \\
 \dl  \: & \ds  \xra{}  \ds \Smash \ds.
\end{array}
\]
The map $\sg$ is obtained by applying $D$ to the fold map $S^0\Wedge
S^0\xra{}S^0$, or by applying $B$ to the coproduct functor
$\amalg\:\CC^2\xra{}\CC$.  The map $\mu$ is obtained by applying $B$
to the product functor $\tm\:\CC^2\xra{}\CC$.  The map $\tht$ only
exists stably.  It is obtained by applying $D$ to the stable pinch map
$S^0\xra{}S^0\Wedge S^0$, or by taking the transfer associated to the
covering map $\sg$.  Finally, $\dl$ is just the diagonal map of the
space $\BC$.  The interaction between these maps is studied in detail
in~\cite{sttu:rme}.

In order to transfer our results about $\ev_*\bs{k}$ to $\ds$, we need
the following lemma.
\begin{lemma}
 Let $M$ be an $E$-module spectrum such that the sequence
 $\{u_0,\ldots,u_{n-1}\}$ is regular on $\pi_*M$.  Then
 $\pi_*\Lhat M=\pi_*(M)^\wedge_{I_n}$.
\end{lemma}
\begin{proof}
 By~\cite[Theorem 2.1 and Corollary 2.2]{ho:blf}, we see that
 $\Lhat M$ is the homotopy inverse limit of spectra $S/I\Smash L_nM$,
 as $S/I$ runs over a suitable tower of generalised Moore spectra of
 type $n$.  As $M$ is a module over the $E(n)$-local spectrum $E$, we
 see that $L_nM=M$.  If $I=(v_0^{a_0},\ldots,v_{n-1}^{a_{n-1}})$, then
 regularity implies that $\pi_*(S/I\Smash M)=\pi_*(M)/I$.  The Milnor
 sequence now implies that $\pi_*\Lhat M=\pi_*(M)^\wedge_{I_n}$, as
 claimed.   
\end{proof}
In particular, the lemma applies when $\pi_*(M)$ is a free module over
$E_*$.  To understand what it says in that case, observe that the
completion of $\bigoplus_k E^0$ is the set of sequences
$\un{a}\in\prod_k E^0$ such that $a_k\xra{}0$ in the $I$-adic topology
as $k\xra{}\infty$.

It is now easy to deduce the following.
\begin{proposition}\label{prop-kunneth}
 $\eds$ and $E^0(\ds\Smash\ds)$ are products of countably many copies
 of $E^0$, and $E^0(\ds\Smash\ds)=\eds\hote\eds$.  Also, $\evds$
 and $\ev(\ds\Smash\ds)$ are completions of countably generated free
 modules over $E^0$, and $\ev(\ds\Smash\ds)=\evds\hote\evds$.
 Moreover, we have $\eds=\Hom_E(\evds,E^0)$. \qed
\end{proposition}

We therefore get (co)products as follows:
\[
\begin{array}{rll}
 *     \: &\evds \hote\evds \xra{}  \evds &\text{ induced by } \sg\\
 \psit \: &\evds \xra{} \evds \hote \evds &\text{ induced by } \tht\\
 \psib \: &\evds \xra{} \evds \hote \evds &\text{ induced by } \dl
\end{array}
\]
\[
\begin{array}{rll}
 \psis \: &\eds \xra{}  \eds \hote\eds &\text{ induced by } \sg\\
 \times\: &\eds \hote \eds \xra{} \eds &\text{ induced by } \tht\\
 \blob \: &\eds \hote \eds \xra{} \eds &\text{ induced by } \dl
\end{array}
\]

\begin{theorem}
 $\eds$ is a Hopf ring with products $\times$ and $\blob$, and
 coproduct $\psis$.
\end{theorem}
\begin{proof}
 We need various diagrams of cohomology groups to commute.  Each one
 is obtained by applying $E^0$ to an obvious diagram of spectra.
 These diagrams are proved to commute in~\cite[Theorem 3.2]{sttu:rme}.
\end{proof}

\section{Morava homology of {$\ds$}}\label{sec-kvds}

Kashiwabara proved in~\cite{ka:bpo} that $\kvds$ is polynomial under
$*$.  It seems worthwhile to give the following simplification of part
of his argument.
\begin{proposition}[Kashiwabara]\label{prop-kvds-poly}
 $\kvds$ is a polynomial algebra under $*$.
\end{proposition}
\begin{proof}
 We can make $\kvds$ into a graded vector space by putting $K_0\bs{m}$
 in degree $m$.  Because $\sg$ sends $\bs{k}\times\bs{l}$ to
 $\bs{k+l}$ and $\tht$ sends $\bs{m}$ to
 $\bigvee_{m=k+l}(\bs{k}\tm\bs{l})_+$, we see that $\kvds$ becomes a
 graded Hopf algebra over $\Fq$ using $*$ and $\psit$.  It is clearly
 connected, in the sense that the degree zero part is just $\Fq$.  It
 follows from Theorem~\ref{thm-EBSk} that $\kvds$ has finite type.  It
 follows in turn from a theorem of Borel that $\kvds$ is isomorphic as
 a graded ring to a tensor product of polynomial and truncated
 polynomial algebras.  The proposition will therefore follow if we
 show that $\kvds$ has no nontrivial nilpotent elements.

 To see this, recall that we have a group-completion map
 $\ds\xra{}QS^0$ and a unit map $QS^0=\Oi S^0\xra{}\oibp$.  The
 composite map $\ds\xra{}\oibp$ induces a ring map
 $\kvds\xra{}K_0\oibp$.  It is shown in~\cite{ka:bpo} that the dual
 map $K^0\oibp\xra{}\kds$ is surjective, so that $\kvds\xra{}K_0\oibp$
 is injective.  It is also known from work of Wilson that $K_0\oibp$
 is a tensor product of polynomial rings and Laurent series rings
 (use~\cite[Theorem 3.3]{wi:osbi} and a degenerate Atiyah-Hirzebruch
 spectral sequence, or see~\cite{rawi:hrc}).  It follows that $\kvds$
 has no nontrivial nilpotent elements, as required.
\end{proof}

\begin{corollary}
 $K_0QS^0$ is the tensor product of a polynomial algebra with the
 group ring 
 \[ \Fq[\Zh]=\Fq[\pi_0 QS^0]\simeq\Fq[u,u^{-1}]. \]
\end{corollary}
\begin{proof}
 It is well-known that $QS^0$ is the group-completion of $CS^0$.  Let
 $Q_0S^0$ be the component of the basepoint in $QS^0$, so that $QS^0$
 is a disjoint union of copies of $Q_0S^0$ indexed by $\pi_0QS^0=\Zh$.
 There is an evident self-map $s\:QS^0\xra{}QS^0$ which sends the
 $k$'th copy to the $(k+1)$'st.  It is not hard to construct a
 self-map $s\:\ds\xra{}\ds$ compatible with this, and to check that
 the induced map from the telescope $s^{-1}\ds$ to $QS^0$ is a weak
 equivalence.  Moreover, for $a\in\kvds$ we have
 $s_*(a)=[1]*a\in\kvds$.  It follows that $K_0QS^0$ is obtained from
 $\kvds$ by inverting the element $[1]$.  Also, $[1]$ is
 indecomposable in $\kvds$, so we can take it as one of our polynomial
 generators.  The claim follows.
\end{proof}

The next corollary is a straightforward consequence of
Propositions~\ref{prop-kunneth} and~\ref{prop-kvds-poly}.
\begin{corollary}
 $\evds$ is the completion of a polynomial algebra over $E^0$, and
 \[ \kvds=\Fq\hote\evds. \qed \]
\end{corollary}

We next prove that $K^0QS^2$ is a formal power series ring, using an
unpublished argument of Kashiwabara.  

\begin{remark}\label{rem-QStwo}
 Here we use the usual product on $K^0QS^2$, coming from the diagonal
 map.  One might think that there were two products on $K^0QS^2$, by
 analogy with the two products on $K^0DS^0$.  Indeed, the two
 projections $S^2\Wedge S^2\xra{}S^2$ give two maps
 $Q(S^2\Wedge S^2)\xra{}QS^2$, and the resulting map
 $Q(S^2\Wedge S^2)\xra{}QS^2\tm QS^2$ is an equivalence.  We can thus
 use the pinch map $S^2\xra{}S^2\Wedge S^2$ to get a map
 $\tht\:QS^2\xra{}QS^2\tm QS^2$, giving a product on $K^0QS^2$.
 However, it is easy to check that $\tht$ is homotopic to the diagonal
 map, so this is the same product as we considered before.  This
 argument does not show that the two products on $K^0DS^0$ are the
 same, because the map $\tht\:DS^0\xra{}DS^0\Smash DS^0$ only exists
 stably.  
\end{remark}

\begin{remark}
 For our purposes, the ring of formal power series over $\Fq$ on a
 countable set of indeterminates $x_0,x_1,x_2,\ldots$ consists of all
 sums $\sum_\al c_\al x^\al$.  Here $\al$ runs over multi-indices
 $\al=(\al_0,\al_1,\ldots)$ such that $\al_i=0$ for $i\gg 0$, and we
 impose no condition at all on the coefficients $c_\al\in\Fq$.  In
 particular, $\sum_i x_i$ is allowed.  
\end{remark}

\begin{proposition}[Kashiwabara]\label{prop-KQS}
 $K^0QS^2$ is a formal power series algebra over $\Fq$.  Moreover, the
 double suspension map induces isomorphisms
 \[ \Ind(K_0QS^0) \simeq \Prim(K_0QS^2) \]
 and
 \[ \Ind(K^0QS^2) \simeq \Prim(K^0QS^0). \]
\end{proposition}
\begin{proof}
 We use the bar spectral sequence.  Let $Y$ be a connected infinite
 loop space, and $X=\Om Y$, so that $Y=BX$.  We then have a spectral
 sequence of Hopf algebras
 \[ E^2_{st} = \Tor^{K_*X}_{s,t}(K_*,K_*) \Rightarrow K_{s+t}(Y) \]
 \[ d^r\: E^r_{st} \xra{} E^r_{s-r,t+r-1}. \]
 Note that $E^r_{st}=0$ where $s<0$; it follows that everything on the
 line $s=1$ survives to $E^\infty$.  The resulting edge map 
 \[ \Ind(K_t X) = \Tor^{K_*X}_{1,t}(K_*,K_*) \Rightarrow K_{1+t}(Y) \]
 is just the homology suspension map.

 We first apply this with $Y=QS^1$, so that $X=QS^0$.  We know that
 $K_*QS^0$ is a polynomial algebra on generators $\{x_i\}$ in degree
 zero, with one generator inverted.  It follows that
 $\Tor^{K_*QS^0}_{st}(K_*,K_*)$ is exterior on generators
 $y_i\in E^2_{10}$.  These are represented by the homology suspensions
 of the generators of $K_0QS^0$.  Because the spectral sequence is
 multiplicative and the generators survive, there can be no
 differentials at all.  Thus $K_*QS^1$ is the exterior algebra on
 primitive classes in odd degree.  Moreover, the suspension map
 induces an isomorphism $\Ind(K_*QS^0)\simeq\Prim(K_*QS^1)$, and the
 evident composite $\Prim(K_*QS^1)\mra K_*QS^1\era\Ind(K_*QS^1)$ is
 also an isomorphism.

 We now repeat the process, and consider the bar spectral sequence
 \[ \Tor^{K_*QS^1}_{st}(K_*,K_*) \Rightarrow 
    K_*(QS^2)\simeq K_*(DS^2). 
 \]
 The $E^2$ term is a divided power algebra on the suspensions of the
 elements $y_i$, which we shall call $z_i$.  The coproduct is just the
 usual one:
 \[ \psi(z_i^{[m]}) = \sum_{m=j+k} z_i^{[j]}\ot z_i^{[k]}. \]
 The differential $d^r$ maps $E^r_{st}$ to $E^r_{s-r,t+r-1}$ and the
 $E^2$ page is concentrated in the lines where $s+t$ is even, so all
 the differentials vanish.  

 It follows that $K_*(QS^2)$ admits a filtration whose associated
 graded is a divided power Hopf algebra.  Dually, we see that
 $K^*(QS^2)$ admits a filtration whose associated graded is a
 polynomial algebra on primitive generators in even degrees (without
 loss of generality, in degree zero).  If we choose elements 
 $z'_i\in K^0(QS^2)$ lifting these generators, we discover in the
 usual way that $K^0(QS^2)\simeq\fps{\Fq}{z'_i}$.  Moreover, the
 suspension map induces an isomorphism
 \[ \Ind(K_*QS^1)\simeq\Prim(K_*QS^2). \] 
 By combining this with information from our first spectral sequence,
 we see that the double suspension induces an isomorphism 
 \[ \Ind(K_*QS^0)\simeq\Prim(K_*QS^2). \]
 Dually, it induces an isomorphism
 \[ \Ind(K^*QS^2)\simeq\Prim(K^*QS^0). \]
\end{proof}

Using Proposition~\ref{prop-kunneth}, we deduce the following
corollary.  
\begin{corollary}\label{cor-EQS-fps}
 $E^0QS^2$ is a formal power series ring over $E^0$.
\end{corollary}

\section{The Snaith splitting}\label{sec-snaith}

In this section, we recall the fundamental theorem of Snaith, which
gives a splitting $\Sgi QX\simeq\Sgi DX$ for connected spaces $X$.  We
will need to prove an additional property of this map (see
Proposition~\ref{prop-susp}) and thus we will need to give some
details of its construction.

We will need to use the usual combinatorial approximation to $QX$,
which we now describe briefly.  First, we define
\[ F_k = F_k\rinf =
   \{\text{injective maps } \; a \: \{1,\ldots,k\} \mra \rinf \}.
\]
It is well-known that this is a contractible free $\Sg_k$-space, in
other words, a model for $E\Sg_k$.  Next, we define
\[ B_k(\rinf;X) = F_k \tm_{\Sg_k} X^k. \]
The disjoint union $\coprod_kB_k(\rinf;X)$ can be thought of as the
set of pairs $(A,x)$, where $A$ is a finite subset of $\rinf$ and
$x\:A\xra{}X$.  We can impose an equivalence relation on this set by
identifying $(A,x)$ with $(B,x|_B)$ if $B\sse A$ and $x$ sends
$A\setminus B$ to the basepoint of $X$.  We define
\[ C(X) = C(\rinf;X) = \left(\coprod_k B_k(\rinf;X)\right)/\sim. \]
\begin{theorem}\label{thm-cx-qx}
 There is a natural map $t\:C(X)\xra{}Q(X)$, which is a weak
 equivalence when $X$ is connected.  On the other hand, when $X=Y_+$
 there is a natural homeomorphism
 \[ C(Y_+) = \coprod_{k\ge 0} B_k(\rinf;Y) =
     \coprod_{k\ge 0} E\Sg_k\tm_{\Sg_k}Y^k = D(Y_+). 
 \]
 In that case, the resulting map $E\Sg_k\tm_{\Sg_k}Y^k\xra{}Q(Y_+)$ is
 adjoint to the composite 
 \[ \Sgi E\Sg_k\tm_{\Sg_k}Y^k_+ \xra{\text{transfer}}
    \Sgi E\Sg_k\tm_{\Sg_{k-1}} Y^k_+ \xra{\text{proj}} 
    \Sgi Y_+.  \qed
 \]
\end{theorem}
\begin{proof}
 For the first sentence, see~\cite[6.3]{ma:gil} (for example).  The
 rest follows easily from the constructions given there.
\end{proof}

We next define 
\[ C_k(X) = \image(B_k(\rinf;X) \xra{} C(X)). \]
This gives a filtration
\[ 0=C_0(X) \sse X=C_1(X) \sse C_2(X) \sse \ldots \sse C(X), \]
and one can easily see that 
\[ C_k(X)/C_{k-1}(X)=E\Sg_{k+}\Smash_{\Sg_k}X^{(k)}=D_k(X). \]

We can now define the Snaith map.  For this, we note that the
polynomial rings $\real[u]$ and $\real[u,v]$ can both be used as
models of $\rinf$.  A point of $a\in C(\real[u];X)$ is a finite set
$A$ of polynomials in $u$, with each polynomial $f\in A$ given a label
$x_f\in X$.  For any subset $B\sse A$, we have a polynomial
\[ g_B = \prod_{f\in B} (v-f) \in \real[u,v]. \]
Clearly, if $B\neq B'$ then $g_B\neq g_{B'}$, by unique
factorisation.  Note also that $(B,x|_B)$ defines a point of
$C_{|B|}(\real[u];X)$; we write $y_B$ for its image in
$D_{|B|}(\real[u];X)$.  By labelling each $g_B$ with the point $y_B$,
we obtain a point $s'(a)=(\{g_B\st B\sse A\},y)\in C(\real[u,v];DX)$.
It can be shown that this gives a well-defined and continuous map 
\[ s'\: C(\real[u];X) \xra{} C(\real[u,v];DX). \]
After composing with the map $t\:C(\real[u,v];DX)\xra{}QDX$ and taking
adjoints, we obtain a map $s\:\Sgi CX\xra{}\Sgi DX=D\Sgi X$.  One
finds that this sends $\Sgi C_kX$ into $\Sgi\bigvee_{j\le k}D_jX$ and
thus that it induces a map $\Sgi(C_kX/C_{k-1}X)\xra{}\Sgi D_kX$, which
is an equivalence.  This is essentially the proof of the following
theorem; for more details, see~\cite{comata:scs}.
\begin{theorem}\label{thm-snaith}
 The Snaith map $s\:\Sgi CX\xra{}\Sgi DX=D\Sgi X$ is a weak
 equivalence.  In the case $X=Y_+$, it can also be obtained by
 applying $\Sgi$ to the homeomorphism $C(Y_+)=D(Y_+)$ mentioned in
 Theorem~\ref{thm-cx-qx}.  \qed
\end{theorem}

There is an obvious equivalence $QS^0\simeq\Om^2 QS^2$, which gives by
adjunction a map $e\:\Sg^2 QS^0\xra{}QS^2$.  We will need to see how
this interacts with the Snaith splitting.  Firstly, we have a map
$e'\:\Sg^2CS^0\xra{}CS^2$, defined as follows: if $z\in S^2$ and 
$A\in CS^0$ is a finite subset of $\rinf$, then $e'(z,A)$ is just the
subset $A$ with each point labelled by $z$.  

Next, recall that $V_k$ is the bundle over $\bs{k}$ corresponding to
the usual representation of $\Sg_k$ on $\cplx^k$.  We write
$\bs{k}^{V_k}$ for the corresponding Thom space.  One can show
directly (see~\cite[Section IX.5]{lemast:esh}) that
\[ D_k(S^2) = \bs{k}^{V_k} = D_k(S^0)^{V_k}. \]
There is an evident vector $v=(1,\ldots,1)\in V_k$ that is fixed under
the action of $\Sg_k$.  Write $L=\cplx v=V_k^{\Sg_k}<V_k$, so that
$\bs{k}^L=\Sg^2\bs{k}$.  The inclusion $L\xra{}V_k$ thus gives a map
$\Dl\:\Sg^2 D_k(S^0)\xra{}D_k(S^2)$ and thus a map
$\Dl\:\Sg^2\ds\xra{}D(S^2)$.  

We will also need a map $\ep_1\:\ds\xra{}S^0$.  This is the counit map
for the coproduct map $\dl\:\ds\xra{}\ds\Smash\ds$.  It sends each
space $\bs{k}$ in $\ds=\bigvee_{k\ge 0}\bsp{k}$ to the non-basepoint
in $S^0$.

\begin{proposition}\label{prop-susp}
 The following diagram commutes.
 \begin{diag}
  \node{\Sg^{\infty+2}QS^0}
  \arrow[3]{e,t}{e}
  \node[3]{\Sgi QS^2} \\
  \node{\Sg^{\infty+2}CS^0} 
  \arrow{n,l}{\Sg^{\infty+2}t}
  \arrow{s,l,=}{\Sg^2 s}
  \arrow[3]{e,t}{e'}
  \node[3]{\Sgi CS^2}
  \arrow{n,lr}{\simeq}{\Sgi t}
  \arrow{s,lr}{\simeq}{s} \\
  \node{\Sg^{\infty+2}DS^0}
  \arrow{e,b}{\Sg^2\tht}
  \node{\Sg^{\infty+2}DS^0\Smash DS^0}
  \arrow{e,b}{1\Smash\ep_1}
  \node{\Sg^{\infty+2}DS^0}
  \arrow{e,b}{\Sgi\Dl}
  \node{\Sgi D(S^2)} 
 \end{diag}
\end{proposition}
\begin{proof}
 It is immediate from the definition of $t$ that the top rectangle
 commutes.  We next observe that the connectivity of $D_m(S^2)$ is at
 least $2m-1$, so that
 $\bigoplus_m\pi_t\Sgi D_m(S^2)=\prod_m\pi_t\Sgi D_m(S^2)$ for all
 $t$.  It follows that the spectrum $D(S^2)$ is the product of the
 spectra $D_m(S^2)$, and of course $\Sg^{\infty+2}\ds$ is the
 coproduct of the spectra $\Sg^{\infty+2}\bsp{k}$.  It is thus enough
 to check that the bottom rectangle commutes after composing with the
 inclusion $\Sg^2\bsp{k}\xra{}\Sg^2\ds$ and the projection
 $D(S^2)=\prod_mD_m(S^2)\xra{}D_m(S^2)$. 

 The two ways around the rectangle give two maps
 $\Sg^{\infty+2}\bsp{k}\xra{}\Sgi D_m(S^2)$.  We shall assume that
 $m\leq k$, leaving it to the reader to see that both maps are zero
 when $m>k$.  We can also take adjoints and identify $QD_m(S^2)$ with
 $C(\real[u,v];D_m(S^2))$ to get two maps of spaces
 $\Sg^2\bsp{k}\xra{}C(\real[u,v];D_m(S^2))$.  We write $F$ for the one
 which comes from $s\circ e'\:\Sg^{\infty+2}\ds\xra{}\Sgi D(S^2)$, and
 $G$ for the one which comes from
 $\Dl\circ(1\Smash\ep_1)\circ\tht\circ t$.
 
 Consider a point $(z,A)\in\Sg^2\bsp{k}$, so $z\in S^2$ and $A$ is a
 finite subset of $\real[u]$ of order $k$.  We see easily from our
 definition of the Snaith map that $F(z,A)\in C(\real[u,v];D_m(S^2))$
 is the collection of polynomials $g_B$ (where $B\sse A$ and $|B|=m$)
 labelled with the points $\Delta(z,B)\in D_m(S^2)$.
 
 Next, recall that $\tht$ is built from the transfer maps
 $\Sgi\bsp{(i+j)}\xra{}\Sgi(\bs{i}\tm\bs{j})_+$.  It follows that we
 have a commutative diagram
 \begin{diag}
  \node{\Sgi\bsp{k}}
  \arrow{s,V} 
  \arrow{e,t}{\text{transfer}}
  \node{\Sgi(\bs{m}\tm\bs{k-m})_+}
  \arrow{e,t}{\text{proj}}
  \node{\Sgi\bsp{m}} \\
  \node{\ds} 
  \arrow{e,b}{\tht}
  \node{\ds\Smash\ds}
  \arrow{e,b}{1\Smash\ep_1}
  \node{\ds}
  \arrow{n,A}
 \end{diag}
 To understand the top composite, we need to recall something about
 transfers.  Let $q\:X\xra{}Y$ be a finite covering map, and let
 $i\:X\xra{}\rinf$ be a map that is injective on each fibre of $q$.
 We then get a map $r\:Y_+\xra{}C(X_+)$ by sending a point $y\in Y$ to
 the set $iq^{-1}\{y\}$, with each point $i(x)\in iq^{-1}\{y\}$
 labelled by $x$.  It is easy to see from the standard constructions
 that $t\circ r\:Y_+\xra{}Q(X_+)$ is adjoint to the transfer map
 $q^!\:\Sgi Y_+\xra{}\Sgi X_+$.

 We now apply this with $X=F_k(\real[u])/(\Sg_m\tm\Sg_{k-m})$ (which
 is a model for $\bs{m}\tm\bs{k-m}$) and $Y=F_k(\real[u])/\Sg_k$
 (which is a model for $\bs{k}$).  We can also think of $Y$ as the set
 of subsets $A\subset\real[u]$ of order $k$, and $X$ as the set of
 pairs $(B,C)$ of disjoint subsets of $\real[u]$ with $|B|=m$ and
 $|C|=k-m$.  The map $q$ just sends $(B,C)$ to $B\cup C$.  For our map
 $i\:X\xra{}\rinf=\real[u,v]$, we take
 $i(B,C)=g_B=\prod_{f\in B}(v-f)$.  The adjoint to the transfer is
 thus the map $Y_+\xra{}C(\real[u,v];X_+)$ that sends $A$ to the set
 of polynomials $g_B$ (for all subsets $B\sse A$ of order $m$), with
 $g_B$ labelled by $(B,C)$.  When we compose with the projection
 $X=\bs{m}\tm\bs{k-m}\xra{}\bs{m}$, we just replace the label
 $(B,C)\in X$ by $B\in\bs{m}$.

 It now follows that adjoint of the double suspension of the top
 composite in the above diagram sends a point $(z,A)\in\Sg^2\bs{k}$ to
 the set of polynomials $g_B$ (for $B\sse A$ with $|B|=k$), labelled
 by $(z,B)$.  To get $G(z,A)$, we just have to apply $\Dl$ to each
 label.  By comparing this with our analysis of $F$, we see that
 $G(z,A)=F(z,A)$ as claimed.
\end{proof}

\section{The Thom isomorphism}\label{sec-thom}

We now discuss the Thom isomorphism.  We saw above that
$D_m(S^2)=\bs{m}^{V_m}$, so we have a Thom isomorphism
$\tilde{E}^0D_m(S^2)=E^0\bs{m}$ and thus $E^0D(S^2)=\eds$.  We need to
show that this respects ring structures in a suitable sense.  Just as
in Remark~\ref{rem-QStwo}, we find that there is only one sensible
product on $E^0D(S^2)$, which comes both from the diagonal map
$D(S^2)\xra{}D(S^2)\tm D(S^2)$ and from the pinch map
$S^2\xra{}S^2\Wedge S^2$.  We give $E^0D(S^0)$ the $\tm$-product,
coming from the stable pinch map $S^0\xra{}S^0\Wedge S^0$.
\begin{lemma}\label{lem-thom-mult}
 With the above product structures, the Thom isomorphism
 $E^0D(S^2)=E^0D(S^0)$ is a ring map.
\end{lemma}
\begin{proof}
 For any spectra $X$ and $Y$, we have a natural map
 $\dl\:D_k(X\Smash Y)\xra{}D_k(X)\Smash D_k(Y)$, which is induced in
 an evident way by the diagonal map $E\Sg_k\xra{}E\Sg_k\tm E\Sg_k$
 when $X$ and $Y$ are suspension spectra.  See~\cite[Section
 I.2]{brmamcst:hir} for discussion of this, and~\cite[Section
 VII.1]{lemast:esh} for proofs.  In particular, there is a natural
 map $\dl\:D_k(\Sg^2X)\xra{}D_k(S^2)\Smash D_k(X)$.  Let
 $u\in E^0D_k(S^2)$ be the Thom class, and define
 $\phi\:E^0D_k(X)\xra{}E^0D_k(\Sg^2X)$ by $\phi(x)=\dl^*(u\ot x)$.
 When $X$ has the form $\Sgi Y_+$ we notice that $D_k(\Sg^2Y_+)$ is
 the Thom space for the pullback of $V_k$ to $D_k(Y_+)$, and that
 $\phi$ is the Thom isomorphism.  The upshot is that this Thom
 isomorphism is natural for stable maps of $\Sgi Y_+$, not merely for
 unstable maps of $Y$.  The lemma follows by considering the stable
 pinch map $S^0\xra{}S^0\Wedge S^0$.
\end{proof}

\begin{corollary}\label{cor-EDS-formal}
 $\eds$ is a formal power series ring under $\tm$.
\end{corollary}
\begin{proof}
 By the lemma, $\eds$ under $\tm$ is isomorphic to $E^0DS^2$ under the
 product coming from the (unstable) pinch map
 $S^2\xra{}S^2\Wedge S^2$.  By the Snaith splitting, $E^0DS^2$ is
 isomorphic to $E^0QS^2$.  By naturality, this converts our product to
 the product on $E^0QS^2$ that comes from the pinch map.  By
 Proposition~\ref{prop-KQS} and Remark~\ref{rem-QStwo}, $E^0QS^2$ is a
 formal power series ring under this product.
\end{proof}

\section{Primitives and indecomposables}\label{sec-prim-ind}

In this section we consider the primitives and indecomposables in the
Hopf ring $\eds$.  We need to interpret the indecomposables in a
completed sense.
\begin{definition}
 We topologise $E^0\bs{k}$, $\ev\bs{k}$ and $\evds$ using the
 $\mxi$-adic topology, where $\mxi=(u_0,\ldots,u_{n-1})$ is the
 maximal ideal in $E^0$.  We give $\eds=\prod_kE^0\bs{k}$ the product
 topology.  Given an augmented topological algebra $A$ over $E^0$, we
 define $\Ind(A)$ to be the quotient of the augmentation ideal by the
 closure of its square.
\end{definition}
\begin{remark}
 With this definition, $\Ind(\fps{E^0}{x_0,x_1,\ldots})$ is a
 countable product of copies of $E^0$, indexed by the variables $x_i$.
\end{remark}

We need to distinguish between two kinds of primitives in $\eds$.
Recall that we have two unit elements $[0]$ and $[1]$ in
$\eds=\prod_kE^0\bs{k}$.  The component of $[0]$ in $E^0\bs{k}$ is $1$
if $k=0$ and $0$ if $k>0$, whereas the component of $[1]$ is $1$ for
all $k$.  We have $[0]\tm a=a=[1]\blob a$ for all $a\in\eds$.  Note
also that $[1]$ is invertible under the $\tm$-product.  Its inverse is
$[-1]$, which is the image of $[1]$ under the map
$D(-1)\:\ds\xra{}\ds$.
\begin{definition}
 We write
 \begin{align*}
  \Prim(\eds)  &= \{x\in\eds\st\psis(x) = x\ot [0]+[0]\ot x\} \\
  \Prim'(\eds) &= \{x\in\eds\st\psis(x) = x\ot [1]+[1]\ot x\}. 
 \end{align*}
\end{definition}
The first of these is the more usual definition for primitives in a
Hopf ring.  It is easy to see from general thoughts about Hopf rings
that the $\blob$-product makes $\Ind(\eds)$ into a ring and
$\Prim(\eds)$ into a module over $\Ind(\eds)$.

On the other hand, one would like to say that the primitives in $\eds$
are dual to $\Ind(\evds)$.  Here $\evds$ is a Hopf algebra with
coproduct $\psis$, and it is natural to use the counit for $\psis$ as
our augmentation on $\evds$.  This counit is $[1]\in\eds$ rather than
$[0]$, and because of this the natural duality is between
$\Prim'(\eds)$ and $\Ind(\evds)$.  This becomes clearer by examining
the explicit formulae in~\cite[Section 6]{sttu:rme}.  As mentioned
there, the map $a\mapsto a\tm [1]$ gives an isomorphism
$\Prim(\eds)\xra{}\Prim'(\eds)$.  

We will need to relate our primitives and indecomposables in terms of
transfers and restrictions, which leads us to our next definitions.
\begin{definition}
 Write $R_m=E^0\bs{m}$.  We shall call the subgroups
 $\bs{i}\tm\bs{j}\le\bs{i+j}$ (with $i,j>0$) {\em partition
   subgroups}.  For each such subgroup $H\le\bs{m}$, we have a
 transfer map $E^0BH\xra{}R_m$; let $I=I_m\leq R_m$ be the sum of the
 images of these maps.  We also have restriction maps
 $R_m\xra{}E^0BH$; let $J=J_m\leq R_m$ be the intersection of the
 kernels of these maps.
\end{definition}

The following two results summarise our conclusions.  They are proved
after Lemma~\ref{lem-IJ-pk}.
\begin{theorem}\label{thm-ind-prim}
 The star-indecomposables and primitives in the Hopf ring $\eds$ are
 as follows:
 \begin{align*}
  \Ind(\eds)  &= \prod_m R_m/I_m = \prod_k R_{p^k}/I_{p^k}      \\
  \Prim(\eds) &= \prod_m J_m = \prod_k J_{p^k}
 \end{align*}
 Moreover, $\Prim(\eds)$ is a free module over $\Ind(\eds)$ on one
 generator $c$, whose component in $R_m=E(\bsp{m})$ is the Euler class
 $c_m=e(V_m-1)$. 
\end{theorem}

\begin{theorem}\label{thm-i-j}
 Suppose that $m=p^k$.  Then $I_m=\ann(c_m)=\ann(J_m)$, and this is a
 retract of $R_m$ as an $E^0$-module.  Moreover, $J_m$ is a free
 module on one generator $c_m$ over $R_m/I_m$ (and thus is generated
 as an ideal by $c_m$).  The ideal $J_m$ is also a retract of $R_m$ as
 an $E^0$-module.  The rank of $R_m/I_m$ over $E^0$ is $\db(k)$ (see
 Definition~\ref{defn-Lambda}).
\end{theorem}

In the light of the above, it is natural to make the following
definition.
\begin{definition}
 $\ov{R}_k=R_{p^k}/I_{p^k}$.
\end{definition}

We now prove a series of lemmas leading to the above theorems.
\begin{lemma}\label{lem-IJ-ideals}
 $I_m$ and $J_m$ are ideals in $R_m$, and $I_mJ_m=0$.
\end{lemma}
\begin{proof}
 This follows from the fact that $\res$ is a ring map, and the formula
 $x\tr(y)=\tr(\res(x)y)$.
\end{proof}

\begin{lemma}\label{lem-c-J}
 The Euler class $c_m=e(V_m-1)$ lies in $J_m$.
\end{lemma}
\begin{proof}
 Suppose that $m=i+j$ with $i,j>0$.  Then the restriction of $V_m$ to
 $\Sg_i\tm\Sg_j$ is isomorphic to $V_i\oplus V_j$, which has two
 trivial summands.  Thus, the restriction of $V_m-1$ has a trivial
 summands and thus $e(V_m-1)$ maps to zero.
\end{proof}

\begin{lemma}\label{lem-IJ-triv}
 If $m$ does not have the form $m=p^k$ then $I_m=R_m$ and $J_m=0$. 
\end{lemma}
\begin{proof}
 In this case it is well-known that
 $(a+b)^m\not\cong a^m+b^m\pmod{p}$, and thus that there exist $i,j>0$
 with $i+j=m$ and $d=|\Sg_m/H|$ not divisible by $p$, where
 $H=\Sg_i\tm\Sg_j$.  The rest of the argument is well-known: we note
 that the composite
 \[ E^0\bs{m} \xra{\res} E^0BH \xra{\tr} E^0\bs{m} \]
 is just multiplication by $\tr(1)$.  The double coset formula shows
 that $\ep(\tr(1))=d\in E^0$, and this is a unit.  Because $E^0\bs{m}$
 is a local ring, we conclude that $\tr(1)$ is a unit.  This implies
 that $\tr\:E^0BH\xra{}E^0\bs{m}$ is surjective, and
 $\res\:E^0\bs{m}\xra{}E^0BH$ is injective.
\end{proof}

\begin{lemma}\label{lem-IJ-pk}
 If $m=p^k$ then 
 \begin{align*}
  I_m &= \image(\tr\:E^0B\Sg_{p^{k-1}}^p\xra{}E^0\bs{p^k}) \\
  J_m &= \ker(\res\:E^0\bs{p^k}\xra{}E^0B\Sg_{p^{k-1}}^p).
 \end{align*}
\end{lemma}
\begin{proof}
 Write $H=\Sg_{p^{k-1}}^p$ and $G=\Sg_{p^k}$.  Consider a partition
 subgroup $L=\Sg_i\tm\Sg_j$, where $i+j=p^k$ and $i,j>0$.  Using the
 argument of Lemma~\ref{lem-IJ-triv}, it will be enough to show that
 $|L/L\cap H|$ is not divisible by $p$, or equivalently that $H$
 contains a Sylow $p$-subgroup of $L$.  There are various ways to see
 this.  The most direct is to use the following well-known formulae
 for the $p$-adic valuations of factorials and binomials:
 \begin{align*}
  v_p(p^k!)            &= (p^k-1)/(p-1)                  \\
  v_p\bsm p^k \\ i\esm &= k - v_p(i) \qquad \text{ if } 0<i<p^k.
 \end{align*}
 If $i=p^{k-1}r+s$ with $0\leq s<p^{k-1}$ then 
 \[ H\cap L=\Sg_{p^{k-1}}^{p-1} \tm \Sg_s \tm \Sg_{p^{k-1}-s} \]
 and we calculate directly that 
 \[ v_p|H\cap L| = (p^k-1)/(p-1) - k + v_p(s) = v_p|L| \]
 as required.
\end{proof}

\begin{proof}[Proof of Theorem~\ref{thm-ind-prim}]
 The augmentation for the Hopf ring $\eds$ is the counit for the
 coproduct $\psis$, which is just the restriction map 
 \[ \eds \xra{} E(D_0S) = E(\bs{0}) = E. \]
 Thus, the augmentation ideal is $\prod_{m>0}E(\bs{m})$.  The star
 product is the product $\times$, which is derived from the map
 $\tht=D(\Delta)\:\ds\xra{}\ds\wedge\ds$.  As mentioned in
 section~\ref{sec-epf}, this sends $\Sgi\bs{m+}$ to
 $\bigvee_{m=i+j}\Sgi(\bs{i}\tm\bs{j})_+$ by the sum of the transfer
 maps.  It follows easily (using Lemma~\ref{lem-IJ-triv}) that
 \[ \Ind(\eds) = \prod_m R_m/I_m = \prod_k R_{p^k}/I_{p^k}. \]
 Similarly, the unit for the star product is the element $[0]$, whose
 component in $E(\bs{m})$ is $1$ if $m=0$ and $0$ otherwise.  The
 coproduct $\psis$ is derived from the map $\sg\:\ds\Smash\ds\xra{}\ds$,
 which sends $\bs{i}\tm\bs{j}$ to $\bs{i+j}$ by the usual map, which
 induces the restriction map in $E$-cohomology.  It follows that an
 element $u\in\eds$ is primitive (i.e. $\psis(u)=[0]\ot u+u\ot [0]$)
 if and only if the component $u_m$ lies in $J_m$ for all $m$.  In
 other words (using Lemma~\ref{lem-IJ-triv} again),
 \[ \Prim(\eds) = \prod_m J_m = \prod_k J_{p^k}. \]
 
 The circle product in $\eds$ is just the ordinary product $\blob$
 induced by the diagonal map of the space $\ds$.  As $I_m$ is an ideal
 in $R_m$, we see that the circle product induces a ring structure on
 $\Ind(\eds)$.  We observed above that $I_mJ_m=0$, which implies that
 $\Prim(\eds)$ is a module over $\Ind(\eds)$.  (Both of these facts
 hold in an arbitrary Hopf ring, and can be proved by elementary
 manipulations).  

 This proves all of Theorem~\ref{thm-ind-prim} except for the last
 sentence.  This is a trivial consequence of Theorem~\ref{thm-i-j},
 which we prove next.
\end{proof}

\begin{proof}[Proof of Theorem~\ref{thm-i-j}]
 Firstly, we know from Corollary~\ref{cor-EDS-formal} that $\eds$ is a
 formal power series ring over $E^0$.  It follows that the ideal
 $R_+^2$ of decomposables is a retract of the augmentation ideal $R_+$
 as an $E^0$-module.  We also know that both ideals are homogeneous
 with respect to the splitting $\eds=\prod_mR_m$.  By choosing
 homogeneous bases we conclude that $I_m=R_+^2\cap R_m$ is a retract
 of $R_m$, as claimed.

 Similarly, we know that $\evds$ is a completed polynomial ring; it
 follows that $\Ind(\evds)$ is a retract of $\evds$.  Dually, we see
 that $\Prim'(\eds)$ is a retract of $\eds$.  By applying the map
 $x\mapsto[-1]\tm x$ (which is an automorphism of $\eds$), we conclude
 that $\Prim(\eds)$ is a retract of $\eds$.  Moreover, $\Prim(\eds)$
 is homogeneous with respect to the splitting $\eds=\prod_mR_m$
 (although $\Prim'(\eds)$ is not).  It follows as in the previous
 paragraph that $J_m$ is a retract of $R_m$.

 In~\cite[Section 6]{sttu:rme} it is shown that there is a natural
 isomorphism $L\ote R_m/I_m\simeq F(\latt_m,L)$, where $L$ is a
 certain extension ring of $E^0$, $\latt_m$ is the set of lattices of
 index $m$ in $\Zp^n$, and $F(\latt_m,L)$ is the ring of functions
 from the finite set $\latt_m$ to $L$ (with pointwise operations).  In
 particular, we see that the rank of $R_m/I_m$ over $E^0$ is the
 number of lattices of index $m$ in $\Zp^n$, which is one definition
 of $\db(k)$ (see Definition~\ref{defn-Lambda}).  It is also shown
 in~\cite{sttu:rme} that $L\ote J_m$ has rank $\db(k)$, and we
 conclude in the same way that $J_m$ has rank $\db(k)$ over $E^0$.

 Next, we show that $\Prim(\eds)$ is freely generated by $c$ over
 $\Ind(\eds)$.  We apply the functor $E^0(-)$ to the diagram in
 proposition~\ref{prop-susp} and use various Thom isomorphisms to
 obtain the following diagram.
 \begin{diag}
  \node{E(QS^2)} \arrow[2]{e,t}{\sg^{*2}} \arrow{s}
  \node[2]{E(QS^0)} \arrow{s} \\
  \node{E(DS^2)} \arrow{e,b}{\blob c}
  \node{E(DS^0)} \arrow{e,b}{\tm [1]}
  \node{E(DS^0)}
 \end{diag}
 (The first map in the bottom row is induced by the inclusions
 $\Sg^2\bs{m}=\bs{m}^\cplx\xra{\Dl}\bs{m}^{V_m}$ and the Thom
 isomorphisms $E(\bs{m}^\cplx)\simeq E(\bs{m})\simeq E(\bs{m}^{V_m})$.
 It is standard that the resulting map is just multiplication by the
 Euler class of $V_m-1$, which is just $c_m$.  This shows that the
 first map on the bottom row is just $x\mapsto x\blob c$, as marked in
 the diagram.)
 
 The right hand vertical map is an isomorphism (by the Snaith
 splitting theorem).  The top arrow is the double suspension map, so
 the image is contained in $\Prim'(\eds)$.  Consider the analogous map
 in $K$-cohomology.  Proposition~\ref{prop-KQS} shows that the image
 of this map is precisely $\Fq\hote\Prim'(\eds)=\Prim'(\kds)$.  As all
 the modules in question are products of free modules, we conclude
 that the double suspension $E(QS^0)\xra{}\Prim'(\eds)$ is surjective.

 Recall that $x\mapsto x\tm [1]$ is an $E^0$-module automorphism of
 $\eds$ which carries $\Prim(\eds)$ isomorphically to $\Prim'(\eds)$.
 It follows from this and the discussion above that the map
 $x\mapsto x\blob c$ is an epimorphism $\Ind(\eds)\xra{}\Prim(\eds)$.
 We therefore have an epimorphism $R_m/I_m\xra{c}J_m$.  Both source
 and target are finitely generated projective (and hence free) modules
 over $E^0$, and we have seen that they have the same rank.  It
 follows that our map $R_m/I_m\xra{c_m}J_m$ is an isomorphism.  This
 means in particular that $\ann(c_m)=I_m$ as claimed.
\end{proof}

\section{Subgroups of formal groups}\label{sec-subgp}

In this section, we identify the functor represented by the ring
$\ov{R}_k$.  We shall assume the results and terminology
of~\cite{st:fsf}, and the discussion of the $E$-cohomology of Abelian
groups in~\cite[Section 4]{grst:vlc}.  We shall write $X=\spf(E^0)$
and $\GG=\spf(E^0\cpi)$, so that $\GG$ is a formal group over $X$.
The special fibre of $X$ is $X_0=\spf(\Fq)$, and the restriction of
$\GG$ over $X_0$ is $\GG_0=\spf(K^0\cpi)$.  Moreover, $\GG$ is the
universal deformation of $\GG_0$.

Now let $Z$ be a space, and $W$ a complex vector bundle over $Z$.  We
can define the associated projective bundle
\[ P(W) = \{(z,L) \st z\in Z \;,\; L \le W_z \;,\; \dim(L)=1 \}.\]
There is a tautological line bundle $L(W)$ over $P(W)$, defined by
\[ L(W)_{(z,L)} = L. \]
This is classified by a map $P(W)\xra{}\cpi$.  We thus obtain a map
$P(W)\xra{}\cpi\tm Z$ and hence $\spf(E^0P(W))\xra{}\GG\tX\spf(E^0Z)$.
If $x$ is the Euler class of $L(W)$ and $a_i$ is the $i$'th Chern
class of $W$, then a well-known lemma says that
\[ E^* P(W) = \fps{E^*Z}{x} / \sum_{k=0}^m a_{m-k} x^k. \]
(This is clear when $W$ is trivialisable, and thus by a Mayer-Vietoris
argument when $Z$ is a finite union of reasonable subsets over which
$W$ is trivialisable, and thus for arbitrary $Z$ by a limit argument.
See~\cite[Theorems 7.4 and 7.6]{cofl:rck} or~\cite{qu:epc} for more
discussion.)  This shows that $\BD(W)=\spf(E^0P(W))$ is a divisor of
degree $m$ on $\GG$ over $\spf(E^0Z)$.

We can apply this process to get a divisor $\BD(V_k)$ over
$\spf(E^0\bs{k})$.  We write $Y_k=\spf(\ov{R}_k)$, which is a closed
subscheme of $\spf(E^0\bs{p^k})$, and we let $H_k$ denote the
restriction of $\BD(V_{p^k})$ to $Y_k$.

Recall that there is a scheme $\Div_{p^k}(\GG)$ over $X$ that
classifies divisors of degree $p^k$ on $\GG$ over $X$.  (It is an
enlightening exercise to identify $\Div_{p^k}(\GG)$ with
$\spf(E^0BU(p^k))$.)  We thus have a map
$\spf(E^0\bs{p^k})\xra{}\Div_{p^k}(\GG)$, which classifies
$\BD(V_{p^k})$.  We also recall from~\cite[Theorem 42]{st:fsf} that
there is a closed subscheme $\Sub_k(\GG)\sse\Div_{p^k}(\GG)$ that
classifies subgroup divisors.  We define
\[ Z_k=\spf(E^0\bs{p^k})\tm_{\Div_{p^k}(\GG)}\Sub_k(\GG).  \]
This is the largest closed subscheme of $\spf(E^0\bs{p^k})$ over which
$\BD(V_{p^k})$ is a subgroup divisor.

\begin{proposition}\label{prop-subgp}
 We have $Y_k\sse Z_k$, so that $H_k$ is a subgroup divisor of $\GG$
 over $Y_k$.
\end{proposition}

\begin{proof}
 Let $A$ be an Abelian $p$-subgroup of $\Sg_{p^k}$ that acts
 transitively on $\{1,\ldots,p^k\}$.  It is not hard to see that the
 restriction of $V_{p^k}$ to $A$ is the regular representation.  This
 can also be written as $\bigoplus_{L\in A^*}L$, the direct sum of all
 one-dimensional representations of $A$.  Recall
 from~\cite[Proposition 4.12]{grst:vlc} that $\spf(E^0BA)$ can be
 naturally identified with $\Hom(A^*,\GG)$.  Over this scheme we have
 a tautological map $\phi\:A^*\xra{}\Gm(\Hom(A^*,\GG),\GG)$, and the
 identification is such that $\BD(L)=[\phi(L)]$ for each $L\in A^*$.  It
 follows that the restriction of $\BD(V_{p^k})$ to $\Hom(A^*,\GG)$ is
 the divisor $\sum_{L\in A^*}[\phi(L)]$.

 If we choose a degree zero complex orientation for $E$, we get Euler
 classes $e(M)\in E^0Z$ for each complex line bundle $M$ over a space
 $Z$.  In particular, a character $L\in A^*$ gives a line bundle over
 $BA$ (which we also call $L$) and thus an Euler class
 $e(L)\in E^0BA$.  We also write $x$ for the Euler class of the
 tautological line bundle over $\cpi$, so $x\in\OO_\GG$ can be thought
 of as a coordinate on $\GG$.  All these identifications work out so
 that $x(\phi(L))=e(L)\in E^0BA=\OO_{\Hom(A^*,\GG)}$.

 In~\cite[Proposition 22]{st:fsf}, we defined a closed subscheme
 $\Level(A^*,\GG)$ of $\Hom(A^*,\GG)$.  It follows
 from~\cite[Proposition 32]{st:fsf} that the restriction of
 $\BD(V_{p^k})$ to $\Level(A^*,\GG)$ is a subgroup divisor.  This means
 that the image of the map
 \[ \Level(A^*,\GG) \mra \Hom(A^*,\GG) = \spf(E^0BA) \xra{}
    \spf(E^0\bs{p^k})
 \]
 is contained in $Z_k$.  It is thus enough to show that the union of
 these images is precisely $Y_k$.  This is essentially clear from the
 rational description of $E^0\bs{p^k}$ given in~\cite{sttu:rme}, but
 some translation is required, so we adopt a slightly different
 approach.  

 Recall from~\cite[Theorem 23]{st:fsf} that $\Level(A^*,\GG)$ is a
 smooth scheme, and thus that $D(A)=\OO_{\Level(A^*,\GG)}$ is an
 integral domain.  Using~\cite[Proposition 26]{st:fsf}, we see that
 when $L\in A^*$ is nontrivial, we have $\phi(L)\neq 0$ as sections of
 $\GG$ over $\Level(A^*,\GG)$, and thus $e(L)=x(\phi(L))\neq 0$ in
 $D(A)$.  It follows that that $c_{p^k}=\prod_{L\neq 1}e(L)$ is not a
 zero-divisor in $D(A)$.  On the other hand, if $A'$ is an Abelian
 $p$-subgroup of $\Sg_{p^k}$ which does not act transitively on
 $\{1,\ldots,p^k\}$, then the restriction of $V_{p^k}-1$ to $A'$ has a
 trivial summand, and thus $c_{p^k}$ maps to zero in $D(A')$.  Next,
 we recall the version of generalised character theory 
 described in~\cite[Appendix A]{grst:vlc}.  It is
 proved there that for any finite group $G$, there is a natural
 isomorphism of rings
 \[ p^{-1}E^0BG = \left(\prod_A p^{-1}D(A)\right)^G,
 \]
 where $A$ runs over all Abelian $p$-subgroups of $G$.  As
 $\ov{R}_k=E^0(\bs{p^k})/\ann(c_{p^k})$ and everything in sight is
 torsion-free, we see that $p^{-1}\ov{R}_k$ is the quotient of
 $p^{-1}E^0\bs{p^k}$ by the annihilator of the image of $c_{p^k}$.
 Using our analysis of the images of $c_{p^k}$ in the rings $D(A)$, we
 conclude that
 \[ p^{-1}\ov{R}_k = \left(\prod_A p^{-1}D(A)\right)^{\Sg_{p^k}},
 \]
 where the product is now over all transitive Abelian $p$-subgroups.
 This implies that for such $A$, the map $E^0\bs{p^k}\xra{}D(A)$
 factors through $\ov{R}_k$, and that the resulting maps
 $\ov{R}_k\xra{}D(A)$ are jointly injective.  This means that $Y_k$ is
 the union of the images of the corresponding schemes
 $\Level(A^*,\GG)$, as required.
\end{proof}

The above proposition gives us a classifying map
$Y_k=\spf(\ov{R}_k)\xra{}\Sub_k(\GG)$, such that the pullback of the
universal subgroup divisor over $\Sub_k(\GG)$ is $H_k$.  Equivalently,
we have a classifying map $\OO_{\Sub_k(\GG)}\xra{}\ov{R}_k$.

\begin{theorem}
 The classifying map $\OO_{\Sub_k(\GG)}\xra{}\ov{R}_k$ is an
 isomorphism.
\end{theorem}

\begin{proof}
 By Theorem~\ref{thm-i-j} above and~\cite[Therem 42]{st:fsf}, we
 know that $\OO_{\Sub_k(\GG)}$ and $\ov{R}_k$ are free modules of the
 same rank over $E^0$.  We shall show that the induced map 
 \[ \OO_{\Sub_k(\GG_0)} = \Fq\ote\OO_{\Sub_k(\GG)}
        \xra{} \Fq\ote\ov{R}_k
 \]
 is injective.  We will conclude that it is an isomorphism by
 dimension count; it follows easily that the original map
 $\OO_{\Sub_k(\GG)}\xra{}\ov{R}_k$ is an isomorphism.

 Recall from~\cite[Proposition 56]{st:fsf} that the socle of
 $\OO_{\Sub_k(\GG_0)}$ is generated by a certain element $a'{}^N$,
 where $N=p+\ldots+p^{n-1}$.  We claim that $a'$ maps to $c_{p^k}$ in
 $\Fq\ote\ov{R}_k$.  To see this, let $K$ be the universal subgroup
 divisor over $\Sub_k(\GG)$.  We have
 $\OO_K=\fps{\OO_{\Sub_k(\GG)}}{x}/f_K(x)$ for a uniquely determined
 monic polynomial $f_K$ of degree $p^k$, called the equation of $K$,
 and $a'$ is defined to be the coefficient of $x$ in $f_K(x)$.  As
 $H_k$ is the pullback of $K$, we see that the image of $f_K$ is the
 equation of $H_k=\BD(V_{p^k})$.  This is just the polynomial
 $\sum_{i=0}^{p^k}a_ix^{p^k-i}$, where $a_i$ is the $i$'th Chern class
 of $V_{p^k}$.  We thus find that the image of $a'$ is
 $a_{p^k-1}=e(V_{p^k}-1)=c_{p^k}$, as claimed.

 We showed in Theorem~\ref{thm-EBSk} that $c_{p^k}^{N+1}\neq 0$ in
 $\Fq\ote E^0\bs{p^k}$, but $c_{p^k}I_{p^k}=0$ so $c_{p^k}^N\neq 0$ in
 $\Fq\ote E^0\bs{p^k}/I_k=\Fq\ote\ov{R}_k$.  Thus, the socle of
 $\OO_{\Sub_k(\GG_0)}$ is mapped injectively.  As any nonzero ideal
 meets the socle, we see that the whole ring is mapped injectively, as
 required.
\end{proof}

\section{A basis for {$\OO_{\Sub_m(\GG)}$}}\label{sec-basis}

We finish by finding a basis for $\ov{R}_m=\OO_{\Sub_m(\GG)}$ over
$E^0$.
\begin{definition}
 Write
 \[ \sg(u,v) = \sum_{i=u}^{v-1} p^i = (p^v - p^u)/(p-1). \]
 Let $a_j$ be the $(p^m-p^j)$'th Chern class of $V_{p^m}$.
 Let $\mu$ and $\nu$ be sequences of the form
 \[ 1  =  \mu_0 < \ldots < \mu_r \le n  \]
 \[ 0 \le \nu_0 < \ldots < \nu_r  =  m. \]
 Write
 \[ b = \prod_{i=0}^{r-1} a_{\nu_i}^{\sg(\mu_{i+1},\mu_i)}.
 \]
 We also define $\rho_j$ for $0\le j<m$ as follows: find the unique
 $i$ such that $\nu_{i-1}\leq j<\nu_i$, and then set
 \[ \rho_j = \begin{cases}
     p^{\mu_i} & \text{ if } \mu_i < n \\
     1         & \text{ if } \mu_i = n
    \end{cases}
 \]
 Now write
 \[ C(\mu,\nu) = \{b a^\al\st\forall j\; 0\leq\al_j<\rho_j\} \]
 and 
 \[ C = \coprod_{\mu,\nu} C(\mu,\nu). \]
\end{definition}

\begin{theorem}\label{thm-basis}
 $C$ is a basis for $\ov{R}_m$ over $E^0$.
\end{theorem}

In order to prove this, we need yet more definitions.
In~\cite[Proposition 49]{st:fsf} we introduced certain quotient rings
$E'(k,l)$ of $\Fq\ote\ov{R}_m$, and we need to find bases for these in
order to make an inductive argument. 
\begin{definition}
 Suppose that $0\le k\le m$ and $0<l\le n$.  Let $\mu$ and $\nu$ be
 sequences of the form 
 \[ l  =  \mu_0 < \ldots < \mu_r \le n  \]
 \[ k \le \nu_0 < \ldots < \nu_r  =  m. \]
 Write
 \[ b = \prod_{i=0}^{r-1} a_{\nu_i}^{\sg(\mu_{i+1},\mu_i)}.
 \]
 We also define $\rho_j$ for $0\le j<m$ as follows.  If $j<k$ we put
 $\rho_j=0$.  If $k\le j<m$ we find the unique
 $i$ such that $\nu_{i-1}\leq j<\nu_i$, and then set
 \[ \rho_j = \begin{cases}
     p^{\mu_i} & \text{ if } \mu_i < n \\
     1         & \text{ if } \mu_i = n.
    \end{cases}
 \]
 Now write
 \[ C_{kl}(\mu,\nu) = \{b a^\al\st\forall j\; 0\leq\al_j<\rho_j\} \]
 and 
 \[ C_{kl} = \coprod_{\mu,\nu} C_{kl}(\mu,\nu). \]
\end{definition}

\begin{proof}[Proof of Theorem~\ref{thm-basis}]
 We use the notation of~\cite[Proposition 49]{st:fsf}.  There we
 introduced certain quotient rings $E'(k,l)$ of $\Fq\ote\ov{R}_m$ and
 elements $v,a\in E'(k,l)$ such that
 \begin{align*}
  E'(0,1)   &= \Fq\ote\ov{R}_m                  \\
  E'(m,l)   &= \Fq                              \\
  E'(k,n)   &= \Fq                              \\
  E'(k,l)/a &= E'(k+1,l)                        \\
  E'(k,l)/v &= E'(k,l+1)                        \\
  v a^{p^l} &= 0 \qquad \text{ if } l<n. 
 \end{align*}
 From the definitions, it is not hard to check that $a$ is (up to
 sign) the image of our Chern class $a_k$ in $E'(k,l)$.  Because of
 the first property, it will suffice to show that $C_{kl}$ is a basis
 for $E'(k,l)$, for all $0\le k\le m$ and $0<l\le n$.  This is easy to
 check when $k=m$ or $l=n$.  Suppose that we know that $C_{k+1,l}$ is
 a basis for $E'(k+1,l)=E'(k,l)/a_k$ and that $C_{k,l+1}$ is a basis
 for $E'(k,l+1)=E'(k,l)/v$.  Using the relation $v a_k^{p^l}=0$, we
 conclude that 
 \[ C'_{kl} = \{a_k^j t\st t\in C_{k+1,l} \text{ and } 0\le j<p^l \}
              \amalg \{a_k^{p^l} t \st t \in C_{k,l+1} \}
 \]
 is a spanning set for $E'(k,l)$.  This gives an upper bound for
 $\dim_{\Fq}E'(k,l)$ and thus by induction a bound for
 $\dim_{\Fq}E'(0,1)=\dim_{E^0}\ov{R}_m$.  This is the conclusion
 of~\cite[Proposition 49]{st:fsf}.  It is proved in~\cite[Corollary
 53]{st:fsf} that this bound is sharp, and it follows easily that the
 bounds used at each stage must be sharp, so that $C'_{kl}$ is a basis
 for $E'(k,l)$.  Thus, we need only check the purely combinatorial
 fact that $C'_{kl}=C_{kl}$.  In fact, one can check that
 \[ \{a_k^j t\st t\in C_{k+1,l}(\mu,\nu) \text{ and } 0\le j<p^l \} =
    C_{kl}(\mu,\nu).
 \]
 This shows that the first piece of $C'_{kl}$ is the disjoint union of
 the sets $C_{kl}(\mu,\nu)$ for which $\nu_0>k$.  Next, consider a
 piece $C_{k,l+1}(\mu,\nu)$ of $C_{k,l+1}$.  If $\nu_0=k$, we define
 \[ \nu' = (\nu_0 < \nu_1 <\ldots <\nu_r = m) \]
 \[ \mu' = (l     < \mu_1 <\ldots <\mu_r \leq n). \]
 If $\nu_0>k$, we define instead
 \[ \nu' = (k <\nu_0 < \nu_1 <\ldots <\nu_r = m) \]
 \[ \mu' = (l <l+1   < \mu_1 <\ldots <\mu_r \leq n). \]
 Either way, one can check that
 \[ \{a_k^{p^l} t\st t\in C_{k,l+1}(\mu,\nu)\} = C_{kl}(\mu',\nu'). \]
 Using this, one can conclude that $C'_{kl}=C_{kl}$ as required.
\end{proof}

%\bibliography{bibdata}

\end{document}